\documentclass[a4paper,12pt,reqno]{amsart}
 
\usepackage{geometry}
\usepackage{indentfirst}
\usepackage{amssymb} 
\usepackage{wrapfig}
\usepackage{lscape}
\usepackage{rotating}
\usepackage{floatrow}
\usepackage{float}
\usepackage{graphicx}
\usepackage{url}
\usepackage[table]{xcolor}
\usepackage{array}
\usepackage{cite}
\usepackage{changepage}
\usepackage{multicol}

\usepackage{amsfonts,amsmath,amssymb,amsthm}
\usepackage{graphicx}
\usepackage[linesnumbered]{algorithm2e}
\usepackage{fixmath}
\usepackage{hyperref}
\usepackage{tikz}
\usepackage{ulem} 
\usepackage{makecell} 

\newtheorem{definition}{Definition}

\let\vec\mathbf
\newcolumntype{+}{!{\vrule width 2pt}}
\graphicspath{ {images/} }
 \newlength\savedwidth

\begin{document}

\pagestyle{plain}

\begin{flushleft}
{\Large
\textbf\newline{Globally optimal dense and sparse spanning trees, and their applications} %
}
\newline

Mustafa Ozen\textsuperscript{1}, 
Goran Lesaja\textsuperscript{2},
Hua Wang\textsuperscript{2}
\\
\bigskip
\textbf{1} New Jersey Institute of Technology, Newark, NJ, 07102, USA
\\
\textbf{2} Georgia Southern University, Statesboro, GA 30460, USA
\\
\bigskip

\end{flushleft}

\section*{Abstract}
Finding spanning trees under various constraints is a classic problem with applications in many fields. Recently, a novel notion of ``dense'' (``sparse'') tree,  and in particular spanning tree (DST and SST respectively), is introduced as the structure that have a large (small) number of subtrees, or small (large) sum of distances between vertices. We show that finding DST and SST reduces to solving the discrete optimization problems. New and efficient approaches to find such spanning trees is achieved by imposing certain conditions on the vertex degrees  which are then used to define an objective function that is minimized over all spanning trees of the graph under consideration. Solving this minimization problem exactly may be prohibitively time consuming for large graphs. Hence, we propose to use genetic algorithm (GA) which is one of well known metaheuristics methods to solve DST and SST approximately. As far as we are aware this is the first time GA has been used in this context. We also demonstrate on a number of applications that GA approach is well suited for these types of problems both in computational efficiency and accuracy of the approximate solution. Furthermore, we improve the efficiency of the proposed method by using Kruskal's algorithm in combination with GA.

The application of our methods to several practical large graphs and networks is presented. Computational results show that they perform faster than previously proposed heuristic methods and produce more accurate solutions. Furthermore, the new feature of the proposed approach is that it can be applied recursively  to sub-trees or spanning trees with additional constraints in order to further investigate the graphical properties of the graph and/or network. The application of this methodology on the gene network of a cancer cell led to isolating key genes in a network that were not obvious from previous studies.

\section*{Introduction}
Seeking the spanning tree of a given graph structure is a classic problem that has numerous applications and variations. For some examples of such study one may see \cite{1.,2.,3.,4.,5.,6.,7.,8.,9.}. In a weighted graph, finding the spanning tree with minimum total weight is known as the minimum spanning tree problem and is probably one of the most extensively studied problems. 

In the case of unweighted graphs, it is of interest to define a formal criterion to distinguish spanning trees (or sub-structures in general) that are more ``compact'' or ``spread out''. One such criteria can be introduced through the topological indices defined as graph invariants. The best known distance-based index is the Wiener index \cite{10.,11.}, defined as the sum of distances between all pairs of vertices. A counting-based index, sometimes called  the subtree index, is defined as the number of subtrees \cite{12.}. It has been observed that in many classes of graphs, the extremal structure that maximizes the Wiener index usually minimizes the subtree index, and vice versa \cite{13.}. Naturally, a tree with many subtrees and small Wiener index is considered as ``dense''  while a tree with few subtrees but large Wiener index is considered as ``sparse''. For example, it is a well known and easily proved fact that among all trees of a given order, the star is the densest and the path is the sparsest tree (see, for instance, \cite{12.}).       

In \cite{6.} an edge-swap heuristics between two spanning trees that finds dense spanning trees through adjacent vertex degrees was proposed based on the theoretical analysis of the problem. 
This topic was further investigated in \cite{9.}, where more general degree conditions were proposed as the criteria to judge the denseness of a spanning tree. The condition is defined as follows.  For a vector (of real numbers) $\vec{j}= \langle j_1,j_2,\ldots, j_i \rangle$,  the condition $C_{\vec{j}}$ is given by
$$ C_{\vec{j}} = C( \langle j_1,j_2,\ldots, j_i \rangle ) = C_{1,j_1} + C_{2,j_2} + \ldots + C_{i, j_i} $$
where 
$$ C_{i,j} = \sum_{d(u,v)=i} \left((deg(u))^j + (deg(v))^j \right)  $$
is the sum of the $j$-th power of pairs of degrees of vertices at distance $i$ apart. 

Throughout the paper, we use the traditional notations, $G$, $T$, $E(G)$, $V(G)$, $d(u,v)$ and $deg(u)$ to denote a graph $G$, a tree $T$, the edge set of $G$, the vertex set of $G$, the distance between two vertices $u$ and $v$, and the degree of a vertex $u$, respectively.

In the case of condition $C_{1,1}$, we simply have the sum of adjacent vertex degrees. A similar expression 
$$C_{\vec{j}} = C({\langle 1,1 \rangle})=\sum_{uv \in E(T)} (deg(u) + deg(v) ) + \sum_{d(u,v)=2} (deg(u) + deg(v) ) $$ is exactly the condition studied in \cite{6.} where the Dense Spanning Tree (DST) problem was solved by maximizing  $C_{\vec{j}}$  over the set of all spanning trees in the given graph.

Alternatively, minimizing a condition such as $C_{\vec{j}}$  can be used to find  a Sparse Spanning Tree (SST) in the given graph. In \cite{9.}, through computational analysis, it was observed that  for
$$ \vec{j} = \langle 4,2,0,0 \rangle  \hbox{ or } \langle 4,2,2,0\rangle  \hbox{ or } \langle 4,2,2,2 \rangle $$
the corresponding objective functions
\begin{equation}\label{obj1}
\sum_{d(u,v)=1} \left((deg(u))^4 + (deg(v))^4 \right) + \sum_{d(u,v)=2} \left((deg(u))^2 + (deg(v))^2 \right) + 2 \cdot \ell_3 + 2 \cdot \ell_4,
\end{equation}
\begin{align} \label{obj2}
& \sum_{d(u,v)=1} \left((deg(u))^4 + (deg(v))^4 \right) + \sum_{d(u,v)=2} \left((deg(u))^2 + (deg(v))^2 \right) \\
& \quad+ \sum_{d(u,v)=3} \left((deg(u))^2 + (deg(v))^2 \right) + 2 \cdot \ell_4 \nonumber
\end{align}
and
\begin{align}\label{obj3}
& \sum_{d(u,v)=1} \left((deg(u))^4 + (deg(v))^4 \right) + \sum_{d(u,v)=2} \left((deg(u))^2 + (deg(v))^2 \right) \\
& \quad + \sum_{d(u,v)=3} \left((deg(u))^2 + (deg(v))^2 \right) + \sum_{d(u,v)=4} \left((deg(u))^2 + (deg(v))^2 \right) \nonumber
\end{align}
appear to be the most effective when solving DST for most graphs. Here $\ell_i$ is the number of pairs of vertices at distance $i$ from each other.

To understand why these choices of conditions stand out, we will briefly introduce the known extremal facts on the Wiener index and the number of subtrees. This will also help us understand how we measure the denseness of a tree. More details of the discussion below can be found in \cite{6.}. 

As mentioned before, the star and path are considered as the densest and the  sparsest trees for good reasons. More interestingly, among tress of given degree sequence, the {\it greedy tree} (defined below) was shown to minimize the Wiener index \cite{14.,15.,16.} and maximize the number of subtrees \cite{17.}. Here the degree sequence is simply the nonincreasing sequence of vertex degrees.

\begin{definition} Given a degree sequence, the greedy tree is constructed
through the following ``greedy" algorithm:

i) Start with a single vertex $v=v_1$ as the root and give $v$ the appropriate number of neighbors so that it has the largest degree;

ii) Label the neighbors of $v$ as $v_2$, $v_3$, $\ldots$, assign to them the largest available degrees such that
$deg(v_2) \geq deg(v_3) \geq \cdots $;

iii) Label the neighbors of $v_2$ (except $v$) as $v_{21}$, $v_{22}$, $\ldots$ such that they take all the largest degrees available and that
$deg(v_{21}) \geq deg(v_{22}) \geq \cdots $, then do the same for $v_3$, $v_4$, $\ldots$;

iv) Repeat (iii) for all the newly labeled vertices, always start with the neighbors of 
the labeled vertex with largest degree whose neighbors are not labeled yet.
\end{definition}

\begin{figure}[ht]

\centering
    \begin{tikzpicture}[scale=0.8]
        \node[fill=black,circle,inner sep=1.5pt] (v) at (10,6) {}; 
        \node[fill=black,circle,inner sep=1.5pt] (v1) at (4,4) {};
        \node[fill=black,circle,inner sep=1.5pt] (v2) at (8,4) {};
        \node[fill=black,circle,inner sep=1.5pt] (v3) at (12,4) {};
        \node[fill=black,circle,inner sep=1.5pt] (v4) at (16,4) {};
        \node[fill=black,circle,inner sep=1.5pt] (v11) at (3,2) {};
        \node[fill=black,circle,inner sep=1.5pt] (v12) at (4,2) {};
        \node[fill=black,circle,inner sep=1.5pt] (v13) at (5,2) {};
        \node[fill=black,circle,inner sep=1.5pt] (v21) at (7,2) {};
        \node[fill=black,circle,inner sep=1.5pt] (v22) at (8,2) {};
        \node[fill=black,circle,inner sep=1.5pt] (v23) at (9,2) {};
        \node[fill=black,circle,inner sep=1.5pt] (v31) at (11,2) {};
        \node[fill=black,circle,inner sep=1.5pt] (v32) at (13,2) {};
        \node[fill=black,circle,inner sep=1.5pt] (v41) at (15,2) {};        
        \node[fill=black,circle,inner sep=1.5pt] (v42) at (17,2) {};         

        \node[fill=black,circle,inner sep=1.5pt] (v111) at (2.7,0) {};
        \node[fill=black,circle,inner sep=1.5pt] (v112) at (3.3,0) {};
        \node[fill=black,circle,inner sep=1.5pt] (v121) at (3.7,0) {};
        \node[fill=black,circle,inner sep=1.5pt] (v122) at (4.3,0) {};
        \node[fill=black,circle,inner sep=1.5pt] (v131) at (4.7,0) {};
        \node[fill=black,circle,inner sep=1.5pt] (v132) at (5.3,0) {};
        \node[fill=black,circle,inner sep=1.5pt] (v211) at (6.7,0) {};
        \node[fill=black,circle,inner sep=1.5pt] (v212) at (7.3,0) {};
        \node[fill=black,circle,inner sep=1.5pt] (v221) at (7.7,0) {};
        \node[fill=black,circle,inner sep=1.5pt] (v222) at (8.3,0) {};
        \node[fill=black,circle,inner sep=1.5pt] (v231) at (9,0) {};        
        \node[fill=black,circle,inner sep=1.5pt] (v311) at (11,0) {};

        \draw (v)--(v1);
        \draw (v)--(v2);
        \draw (v)--(v3);
        \draw (v)--(v4);
        \draw (v1)--(v11);
        \draw (v1)--(v12);
        \draw (v1)--(v13);
        \draw (v2)--(v21);
        \draw (v2)--(v22);
        \draw (v2)--(v23);
        \draw (v3)--(v31);
        \draw (v3)--(v32);
        \draw (v4)--(v41);
        \draw (v4)--(v42);
        \draw (v11)--(v111);
        \draw (v11)--(v112);
        \draw (v12)--(v121);
        \draw (v12)--(v122);
        \draw (v13)--(v131);
        \draw (v13)--(v132);
        \draw (v21)--(v211);
        \draw (v21)--(v212);
        \draw (v22)--(v221);
        \draw (v22)--(v222);
        \draw (v23)--(v231);
        \draw (v31)--(v311);

\node at  (10,6.4) {$v_1$};

\node at  (3.5,4) {$v_2$};
\node at  (7.5,4) {$v_3$};
\node at  (12.5,4) {$v_4$};
\node at  (16.5,4) {$v_5$};

\node at  (2.5,2) {$v_{21}$};
\node at  (3.5,2) {$v_{22}$};
\node at  (4.5,2) {$v_{23}$};
\node at  (6.5,2) {$v_{31}$};
\node at  (7.5,2) {$v_{32}$};
\node at  (8.5,2) {$v_{33}$};
\node at  (10.5,2) {$v_{41}$};
\node at  (13.5,2) {$v_{42}$};
\node at  (14.5,2) {$v_{51}$};
\node at  (17.5,2) {$v_{52}$};

    \end{tikzpicture}

\caption{A greedy tree.}
\label{greedy_pic}
\end{figure}
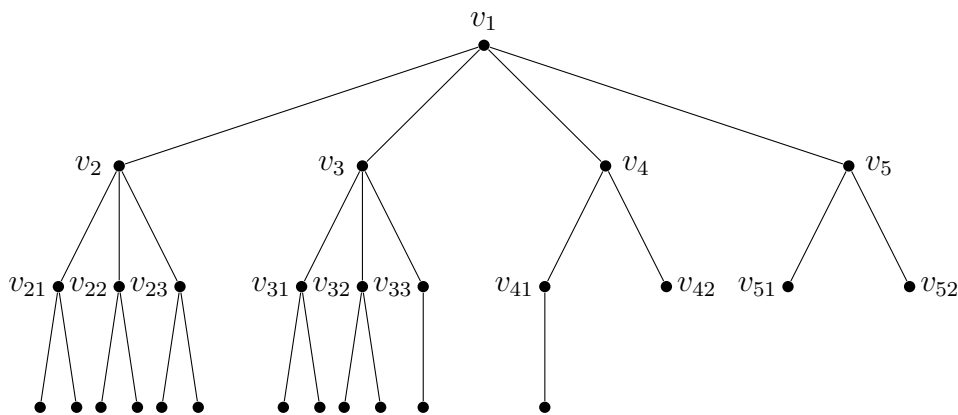

For example, Fig.~\ref{greedy_pic} shows a greedy tree with degree sequence 
$$ ( 4, 4, 4, 3,3,3,3,3,3,3,2,2, 1, \ldots , 1 ) . $$

It is useful to compare the greedy trees of different given degree sequences. In particular, for two nonincreasing sequences $\pi'=(d_1',\cdots,
d'_{n})$ and
 $\pi''=(d''_1, \cdots, d''_{n})$, $\pi''$ is said to {\it majorize}  $\pi'$ if  for $k=1, \cdots, n-1$ we have
 \begin{eqnarray*}
  \sum_{i=0}^{k}d'_i\le\sum_{i=0}^k d''_i \qquad \text{ and } \qquad \sum_{i=0}^{n}d'_i=\sum_{i=0}^{n}d''_i. \end{eqnarray*}
Through the concept of majorization researchers have been able to find the dense structures (with minimal Wiener index or maximal number of subtrees) under various constraints. See \cite{17.} for an example of such discussions. To find dense spanning trees, our edge-swap heuristic starts with a random spanning tree. We then continuously remove a ``bad'' edge and add a ``good'' edge in order to ``improve'' the degree sequence by putting large degrees closer to each other. The corresponding criteria for this procedure are the degree sum conditions discussed earlier. It also makes sense that our optimal conditions given above put more emphasis on the adjacent degree sums than others.

In this paper, we model finding SST and DST as an optimization problem and solve it using the genetic algorithm (GA). Solving minimization problem for SST or maximization problem for DST exactly may be prohibitively time consuming for large graphs. Hence, the well-known heuristic method, the GA is used because it is well suited to solve these types of problems. Genetic Algorithm is a metaheuristic optimization method which attempts to find global maximum  or at least its good approximation. This technique can be applied to unconstrained and constrained problems. In addition, GA is well suited for the problems that are discrete and combinatorial in nature, DST and SST being good examples. There is an extensive literature on GA of which we mention \cite{18.} as a good starting point for further reading. Furthermore, we improve the efficiency of the proposed method by using Kruskal's algorithm in combination with GA. In addition, recursive use of the method to identify certain key nodes (as applied to gene network) reveals new target genes of potential interest to practical medical research. 

The paper is organized as follows. In Section 2, Methodologies, we present the model for finding Dense and Sparse Spanning Trees in the given graph which is then solved using GA approach. We also present a modification of the original model in which the feasible minimum spanning trees are found using Kruskal's algorithm  before the GA is used. Using this approach the feasible set is significantly reduced which results in much faster solution of the problem. In Section 3, Results, we first use a simple objective function
 as an example to illustrate the methodology. We then apply models and methods developed in the  Section 2 to several structures from practical applications and comment on the results. In Section 4, Discussion, we discuss the generalization of this methodology to find various dense or sparse sub-structures or spanning trees under additional  constraints. For a specific application on the gene networks of a cancer cell, we present a recursive application of our algorithm to quickly obtain deeper understanding of the graphical properties of a network. Finally, Section 5, Conclusion, contains concluding remarks and  brief discussions on possible directions for further research.
 
\section{Methodologies}
Given an undirected graph $G$, the goal is to find an acyclic subgraph (i.e. a tree) $T$ which contains all of the nodes in $V(G)$ and optimizes the objective function under consideration. Let $N=|V(G)|$ be the number of nodes in $G$, then the number of edges in a spanning tree $T$ is $|E(T)|=N - 1$. Let all of the edges in $E(G)$ be labeled from 1 to $|E(G)|$, then any subtree $T$ of $G$ can be represented by a vector of edge labels $\vec{h} = \langle h_i \rangle$ ($i = 1,\dots,N-1$). For instance, suppose that we are given the graph  containing  $N = 5$ nodes and 9 edges and the edge labels are as shown in Fig~\ref{fig2}A. Then, $\vec{h} = \langle 1, 3, 4 , 6 \rangle$ represents the tree in Fig~\ref{fig2}B. Using this representation, we propose two models in (\ref{eq:Eq3}) and (\ref{eq:Eq4}) for finding dense or sparse spanning trees:
\begin{figure}[ht]
\centering
\includegraphics[width=110mm]{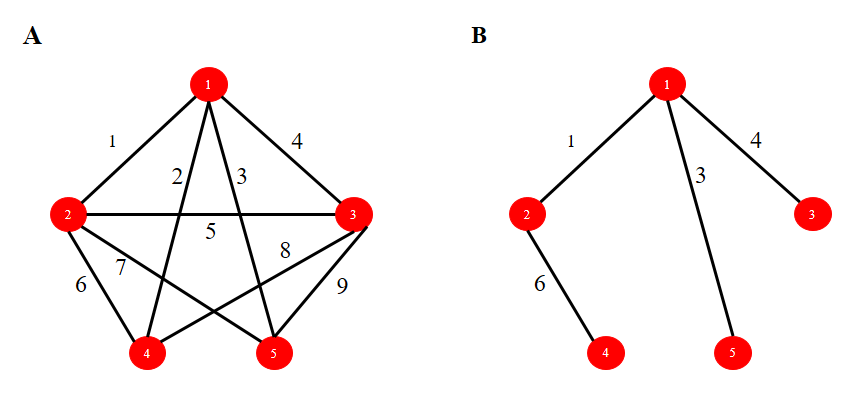}
\caption{{\bf Vector representation of trees.}
(A) Example graph with 5 nodes and 9 edges. (B) Tree represented by $\vec{h} = \langle 1, 3, 4, 6\rangle$.}
\label{fig2}
\end{figure}

\newpage
\subsection{Model 1:}
\begin{itemize}
    \item[]
    {\bf Formulation}

Herein we reduce finding dense or sparse spanning trees problem to the following general form:
\begin{equation}
\label{eq:Eq3}
	\begin{aligned}
		& \underset{\vec{h}}{\text{min}} \text{	Objective Function}(\vec{h}) \\
		\text{subject to}& \\
		& h_i \in \{1,2,\dots,|E(G)|\} \subset \mathbb{Z}^{+}, i = 1,\dots,N-1,  \\
		& h_i \neq h_j, \forall i \neq j, \\
		& \vec{h} \text{ correspond to a connected acyclic subgraph.}
	\end{aligned}
\end{equation}

The objective function in (\ref{eq:Eq3}) could be one of the objective functions in (\ref{obj1})-(\ref{obj3}).
More specifically, we search for an $(N-1)$-dimensional vector $\vec{h}$ with unique integer components $h_i$ such that $\vec{h}$ represents a tree that minimizes the objective function. 

\item[]
{\bf Method:}
Model 1, (\ref{eq:Eq3}), is solved approximately using a well-known metaheuristic method, called Genetic Algorithm (GA). Hence, finding an exact global optimum and corresponding global optimizer(s) is not guaranteed. However, computational experiments show that GA, in most instances, finds  global optimum or a very close approximation of it for these types of problems. The method consists of two phases, ``Pre-optimization'' and ``Optimization'',  described below.

\begin{itemize}
\item[]
{\it Pre-optimization}
\begin{enumerate}
	\item{Create adjacency matrix for $G$. We use the adjacency matrix to keep track of structure information and to plot the optimal trees. To illustrate, $A_1$ and $A_2$ below, where $A_2$ is a variation of  $A_1$ with edge labels as nonzero entries, instead of  1s, are the adjacency matrix of the graph in Fig~\ref{fig2}A. 

\begin{multicols}{2}	
	$$A_1 = \begin{bmatrix}
    	0 & 1 & 1 & 1 & 1 \\
    	1 & 0 & 1 & 1 & 1 \\
    	1 & 1 & 0 & 1 & 1 \\
    	1 & 1 & 1 & 0 & 0 \\
    	1 & 1 & 1 & 0 & 0
	\end{bmatrix} $$ 

	$$ A_2 = \begin{bmatrix}
    	0 & 1 & 4 & 2 & 3 \\
    	1 & 0 & 5 & 6 & 7 \\
    	4 & 5 & 0 & 8 & 9 \\
    	2 & 6 & 8 & 0 & 0 \\
    	3 & 7 & 9 & 0 & 0
	\end{bmatrix}$$
\end{multicols}
	}
	
	\item{Label each edge in $E(G)$ from 1 to $|E(G)|$.}
	\item{Choose an objective function that will be minimized or maximized. Note that any maximization problem can be written as a minimization problem by simply taking the negative of the objective function, and vice versa.}
\end{enumerate}

\item[]
{\it Optimization}
\begin{enumerate}
	\item{Define initial values: We first define initial values for the GA, which are population size, number of variables and variable ranges. Setting population size higher means better chance to obtain globally optimal solution. The number of variables depends on the number of nodes $(N)$ as explained above. The variable ranges depend on the number of edges since the edges are labeled from 1 to $|E(G)|$.}
	\item{Construct the constraints: The first constraint is that the variables are positive integers since they represent the labels of edges. Secondly, the variables should be unique. More specifically, the solution vector cannot contain the same value more than once, i.e. $h_i \neq h_j$ if $i \neq j$. Finally, the vector $\vec{h}$ induces a connected acyclic subgraph. This is controlled by checking the distance between each pair of nodes in the graph represented by $\vec{h}$. If the graph is not connected, then there exists at least one pair of nodes with infinite distance. In such cases, objective function value of the the corresponding $\vec{h}$ is high  and will not survive, that is it will not be included} in the next generation of the GA.
\end{enumerate}
The optimization problem is implemented in MATLAB and solved by using its global optimization toolbox, i.e. the genetic algorithm function ``ga'' \cite{19.}.
\end{itemize}
\end{itemize}

\subsection{Model 2:}
The formulation (\ref{eq:Eq3}) in Model 1 is applicable  not only to relatively small graphs but also to large graphs with hundreds of nodes and edges. However, to obtain a solution for the large graphs, the population size would need to be set very high since the dimension of the problem (depending on the numbers of nodes and edges) increases. Hence, the complexity of the problem increases dramatically in terms of CPU time and number of iterations. 

To address this issue, we  modify  Model 1  by reducing the feasible set  in (\ref{eq:Eq3}) to a set of minimum spanning trees which are constructed using well-known Kruskal's algorithm \cite{20.,21.}. The modified model, Model 2, (\ref{eq:Eq4}) is described below.
\begin{equation}
\label{eq:Eq4}
	\begin{aligned}
		& \underset{\vec{h}}{\text{min}} \text{	Objective Function}(\vec{h}^{*}) \\
		\text{subject to}& \\
		& h_i \in \{1,2,\dots,|E(G)|\} \subset \mathbb{Z}^{+}, i = 1,\dots,N^{*},  \\
		& h_i \neq h_j, \forall i \neq j, \\
		& \vec{h} \text{ containing at least one edge adjacent to each } v \in V(G), \\
		& \vec{h}^{*} = \text{ Kruskal}(\vec{h}).
	\end{aligned}
\end{equation}

In this model $\vec{h}^{*}$ denotes the minimum spanning tree obtained from $\vec{h}$, which resulted from a subset of $E(G)$, and $N^{*}$ is the size of $\vec{h}$. The main idea here is to find a subset $\vec{h}$ of $E(G)$ containing at least one edge adjacent to each node such that the minimum spanning tree $\vec{h}^{*}$ constructed using Kruskal's algorithm from $\vec{h}$ has the optimal objective function value. The advantage of this approach lies in the fact that we narrow down the search space and reduce the cost of checking if the selected $\vec{h}$ is a tree (as it is guaranteed by the Kruskal's algorithm). Thus, we increase the chance of finding a solution with lower population size and consequently reduce the complexity.  

The input of the  Kruskal's algorithm is a graph with weighted edges. However, if  edges do not have weights, then each edge is assigned a weight of 1. Consequently, Model 2,  (\ref{eq:Eq4}) can be applied to weighted and unweighted graphs. 

To find the solution we only need to apply the previous method to the formulation \eqref{eq:Eq4}.

\section{Results}
\subsection{An Introductory Example}
\begin{figure}[ht]
\begin{center}
	\includegraphics[width=125mm]{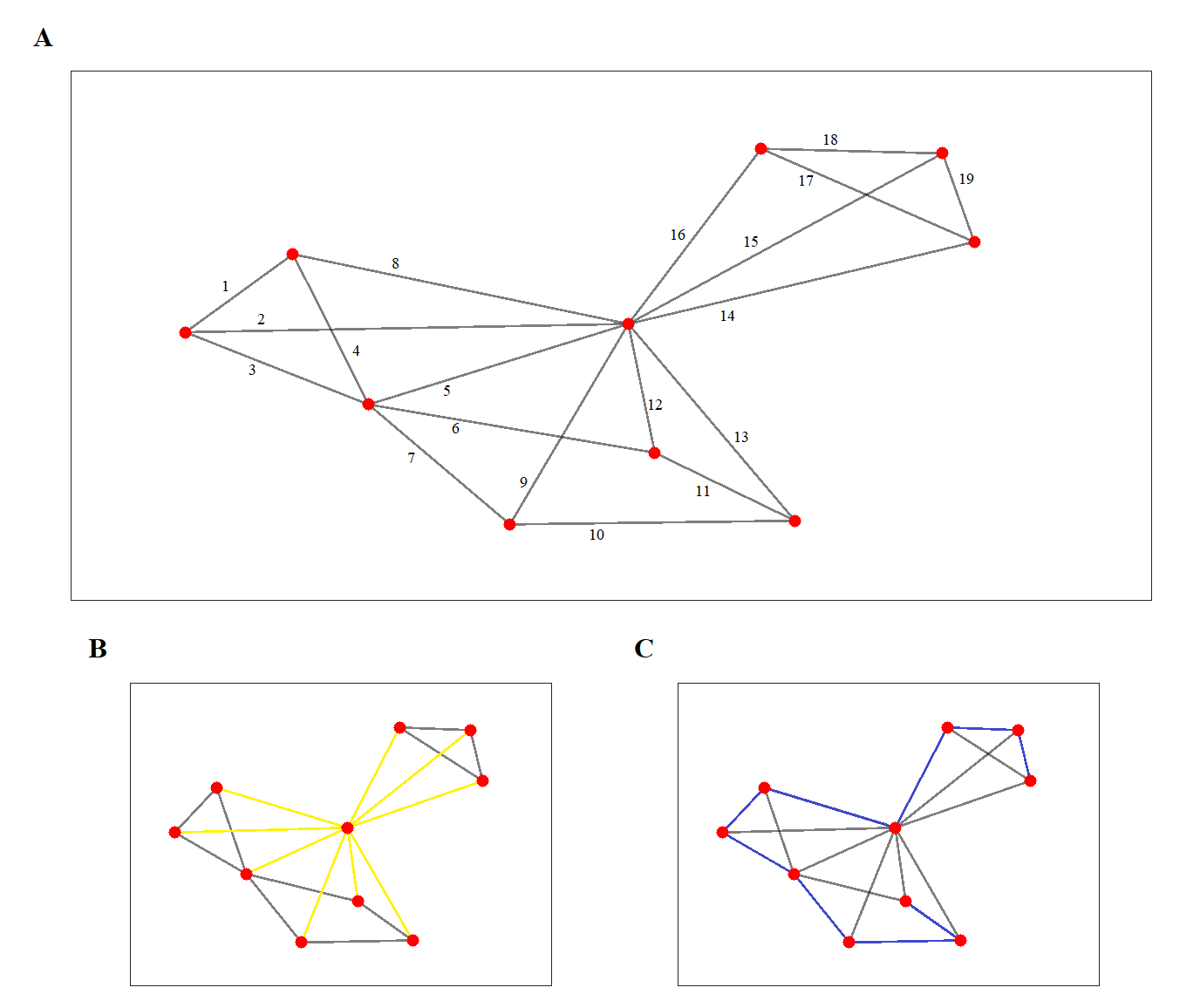}
\end{center}
\caption{{\bf The random graph with 10 nodes and 19 edges.}
(A) Randomly generated undirected graph (B) The ``star'' highlighted in the graph, i.e. $\vec{h} = \langle 2, 5, 8, 9, 12, 13, 14, 15, 16 \rangle$. (C) The ``path'' highlighted in the graph, i.e. $\vec{h} = \langle 1, 3, 7, 8,  10, 11, 16, 18, 19 \rangle$.}
\label{fig3}
\end{figure}

In this subsection we provide a brief justification and an additional argument for the validity of the GA approach described in the previous section. A good way to test our method is to apply it to graphs that contain some obvious dense or sparse spanning trees. Indeed, among trees of the same order it is well known that the star is the densest (maximizing the number of subtrees and minimizing the Wiener index) and the path is the sparsest. Fig~\ref{fig3} shows a randomly generated undirected graph with 10 nodes and 19 edges, containing both the star and the path among its spanning trees. 

The DST of this graph, which is going to be the star, is the global optimal solution of the following Integer Linear Programming (ILP) problem described below. The SST, which is going to be the path, is obtained by finding the minimum of the objective function instead of the maximum. To distinguish ILP model from the proposed models in the previous section, we use the following notation:

Suppose that $G$ represents the undirected graph in Fig~\ref{fig3} with vertex set $V=\{v_j, j=1..n_N\}$ and edge set $E=\{y_i, i=1..n_E\}$ where $n_N$ and $n_E$ are the number of nodes and edges respectively. Let $I_j$ be the corresponding index set of the edges that are connected to vertex $v_j$. Our question can then be formulated as

\begin{equation}
\label{eq:Eq_ILP}
	\begin{aligned}
		& \underset{y}{\text{max}} \ \sum\limits_{k=1}^{n_E} y_k(deg(v_i)+deg(v_j)) \\ 
		\text{subject to}& \\
		& \ deg(v_j) = \sum\limits_{i} y_i, \; i\in I_j, \; \forall j = 1,...,n_N, \\
		& \sum\limits_{i} y_i \geq 1, \; i \in I_j, j = 1,...,n_N, \\
		& \sum\limits_{i} y_i = n_N - 1, \; i = 1,...,n_E, \\
		& 0 \leq y_i \leq 1, \:\: y_i,deg(v_j) \in \mathbb{Z}, \; i = 1,...,n_E, j= 1,...,n_N.
	\end{aligned}
\end{equation}
where $v_i$ and $v_j$ in the objective function are the two end nodes of $y_k$. Solution of this problem is the vector $y$ whose entries, i.e. $y_i$'s, are 1 if and only if an edge should remain in the graph and 0 if the $i^{th}$ edge should be removed to obtain the densest or sparsest tree. Furthermore, each constraints are explained below:

\begin{itemize}
\item The first constraint evaluates the degree of each node where $deg(v_j)$ is the degree of node $v_j$. Value of $deg(v_j)$ depends $y_i$ values $(i \in I_j)$. Recall that $y_i=1$ if the corresponding edge is connected to $v_j$.

\item The second constraint states that every node in the graph will be connected to at least one of the other nodes so that resulting tree will be connected.

\item The third constraint makes sure that the number of edges equals to the number of nodes - 1 so that the result will be a tree. Furthermore, the second and third constraints together makes sure a spanning tree of the graph. 

\item The last constraint makes sure that the variables are integers and $y_j \in \{0,1\}$.  
\end{itemize}

\begin{figure}[ht]
\begin{center}
	\includegraphics[width=120mm]{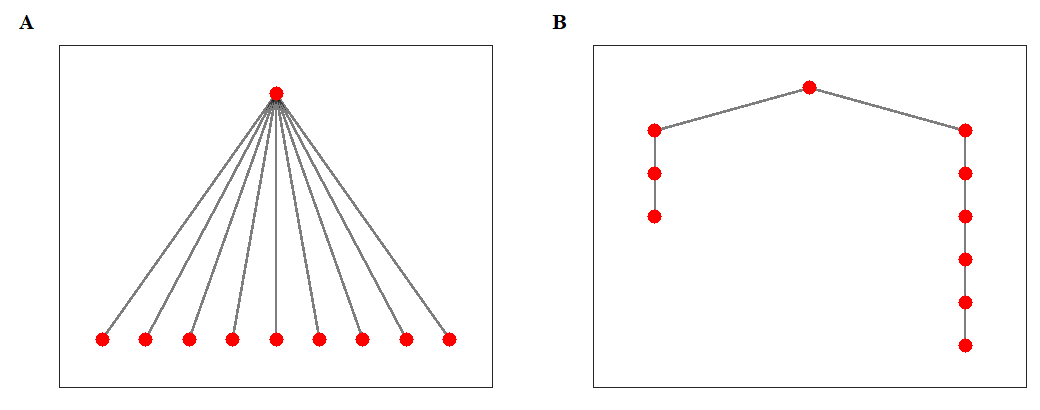}
\end{center}
\caption{{\bf Globally optimal solutions to (\ref{eq:Eq3}) for each objective function.}
(A) The star: optimal solution to minimizing $W(\vec{h})$, $S_{p=1/2}(\vec{h})$ and maximizing $S_{p=2,3}(\vec{h})$, i.e. $\vec{h} = \langle 2, 5, 8, 9, 12, 13, 14, 15, 16 \rangle$. (B) The path: optimal solution to maximizing $W(\vec{h})$, $S_{p=1/2}(\vec{h})$ and minimizing $S_{p=2,3}(\vec{h})$ , i.e. $\vec{h} = \langle 1, 3, 7, 8,  10, 11, 16, 18, 19\rangle$. $W(\vec{h})$ is Wiener index of the tree corresponding to $\vec{h}$.} 
\label{fig4}
\end{figure}

The ILP model (\ref{eq:Eq_ILP}) for the graph in Fig~\ref{fig3}A, is then solved using the well known integer linear  programming solver, IBM ILOG Cplex \cite{25.},  obtaining exact global optimal solution, that is the star  in Fig~\ref{fig4}. The path is obtained by minimization version of the model. 

The same graph in Fig~\ref{fig3}A is then solved, using Model 1 and Model 2 introduced above with the choice of different objective functions. The GA also finds optimal dense and sparse spanning trees (Fig~\ref{fig4}) that are the star (Fig~\ref{fig4}A) and the path (Fig~\ref{fig4}B), respectively. Note, as already mentioned, that minimization version of the models lead to finding SST and maximization versions lead to finding  DST. We use the simple objective function $$ S_p(T) = \sum_{v \in V(T)} (deg(v))^p, $$ 
where $p$ is a fixed positive real parameter, with $p=2$, $p=3$, and $p=1/2$ for GA. In fact, the objective function in (\ref{eq:Eq_ILP}) equals to $S_p(T)$ for $p = 2$. The objective values obtained by GA for different $p$ values are summarized in Table~\ref{table1}.

Although the ILP model (\ref{eq:Eq_ILP}) guarantees the exact global solution, formulating the model and solving it using integer programming methods becomes increasingly more difficult and time consuming  as the size of the graph increases rendering it not practical for large graphs.  On the other hand, models  (\ref{eq:Eq3}) and (\ref{eq:Eq4}) are easier to construct for the relatively large graphs. Furthermore,  GA methodology solves these models efficiently even for large graphs and still provides approximate solutions that are very close to the exact global solutions.

\begin{table}[ht]
\centering
\caption{
{\bf Summary of the results for the graph in the Figure \ref{fig3}.}}
\begin{tabular}{|c|c|c|c|c|}
\hline
 & \multicolumn{2}{|c|}{\bf Dense spanning tree} & \multicolumn{2}{|c|}{\bf Sparse spanning tree}\\ \hline
Objective function & Type & Objective value & Type & Objective value \\ \hline
$S_{p=2}(\vec{h})$ & max & 90 & min & 34 \\ \hline
$S_{p=3}(\vec{h})$ & max & 738 & min & 66 \\ \hline
$S_{p=1/2}(\vec{h})$ & min & 12 & max & 13.3137 \\ \hline
\end{tabular}
\label{table1}
\end{table}
For each $p$, both Model 1 and Model 2 are capable of finding the optimal solution for the given random graph in Fig~\ref{fig3} as summarized in Table~\ref{table1}. However, the methods differ from each other in their efficiency, as discussed in the previous section. More specifically, with Model 1 the solutions were found in about 210 seconds whereas with Model 2 they were found in as quick as 55 seconds for each $p$ (the CPU time may vary depending on the computer specifications). The reason is due to the fact that in Model 2, we do not need to check if the optimizer $\vec{h}$ is a tree or not by employing Kruskal's algorithm. For this specific example Model 2 is almost 4 times faster than Model 1. 

To further examine how well this approach performs, we compare it with our previously proposed heuristic algorithm in \cite{6.}. The heuristic algorithm in \cite{6.} was applied to the US Airports data set, a connected graph of 332 vertices and 2126 edges. For the resulted dense tree, measured by the (small) value of its Wiener index, we obtained the value of 1412038 \cite{6.}. On the other hand, when the GA Model 2 is applied, we found a denser tree with Wiener index as low as 188200. This result indicates that our new method outperforms the previously established heuristic algorithm.

\subsection{Application to practical structures} 

We now apply the proposed methods with objective function $C_{\vec{j}}$ for $\vec{j} = \langle4,2,2,2\rangle$ which is defined in (\ref{obj3}) and proposed in \cite{9.}. Minimization of $-C_{\vec{j}}$ produces dense spanning trees whereas minimization of $C_{\vec{j}}$ results in sparse spanning trees. Most of these structures are relatively large graphs, hence, we use our Model 2 in all of the applications below.

\subsubsection{Gene networks}
\begin{figure}[ht]
\begin{center}
	\includegraphics[width=90mm]{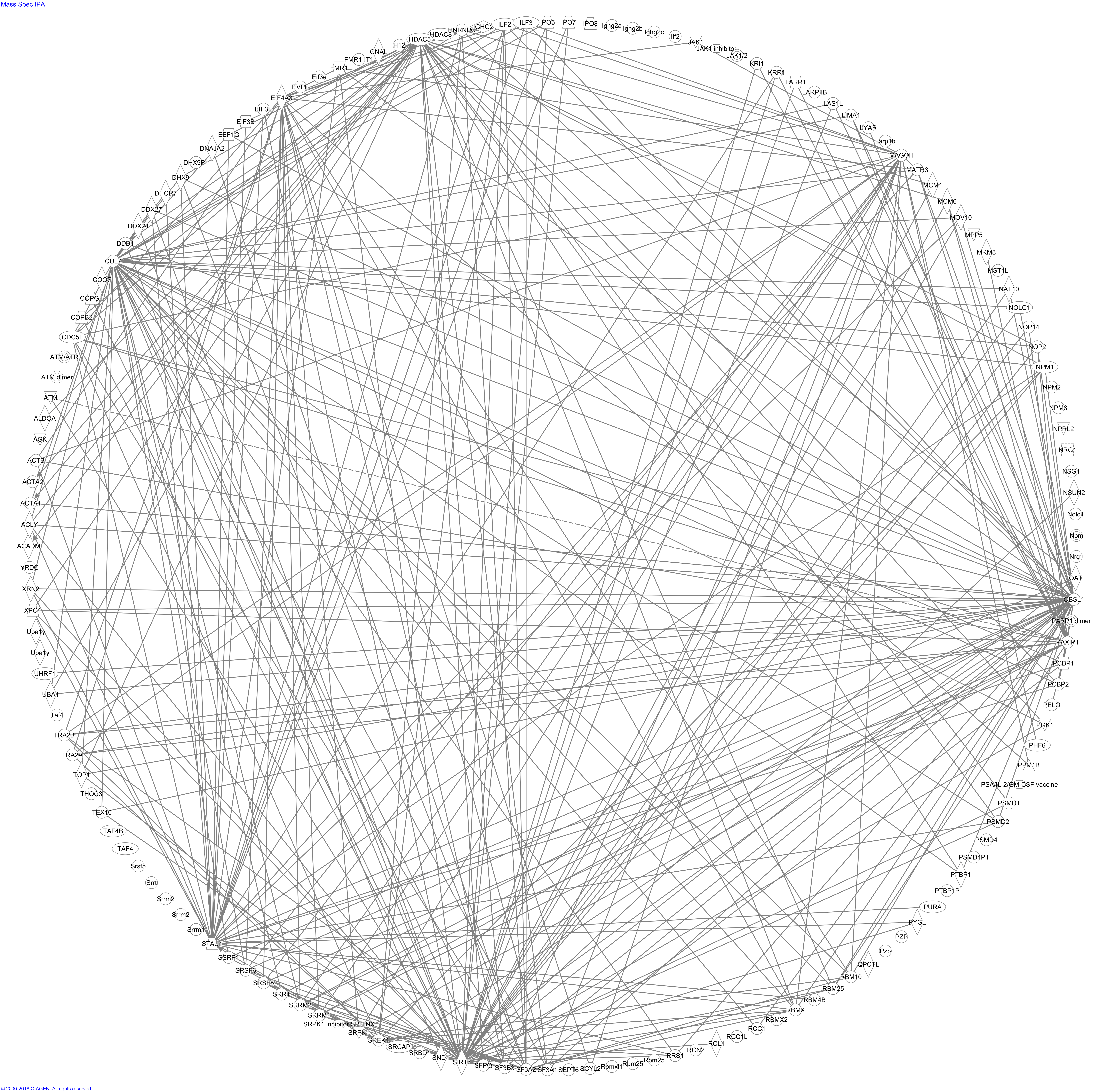}
\end{center}
\caption{{\bf A gene network with 93 genes and 295 interactions.}}
\label{fig5}
\end{figure}

\begin{figure}[ht]
\begin{center}
	\includegraphics[width=100mm]{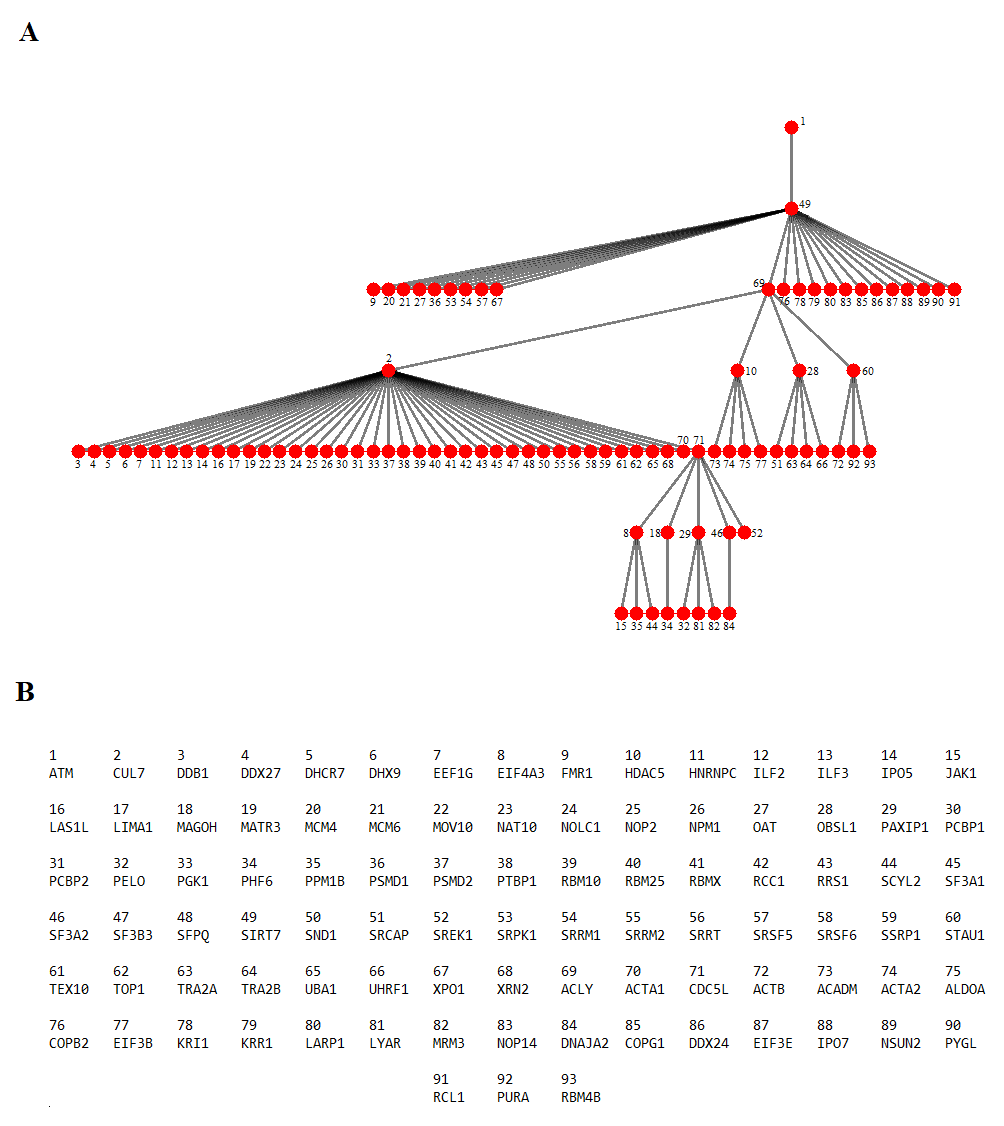}
\end{center}
\caption{{\bf Dense spanning tree of the given gene network.}
(A) The dense spanning tree maximizing $C_{\vec{j}}$ for $\vec{j} = \langle 4,2,2,2 \rangle$. (B) Names and labels of each node.}
\label{fig6}
\end{figure}

The gene network shown in Fig~\ref{fig5} has 93 genes and 295 interactions between them established in the literature. For simplicity we use an undirected edge (of weight 1) to connect any pair of genes (nodes) that have interactions. Using Model 2 with the objective function $C_{\vec{j}}$ for $\vec{j} = \langle 4,2,2,2 \rangle$, we obtain the dense spanning tree in Fig~\ref{fig6}A with objective value 3675370. Gene names and their associated labeles can be found in Fig~\ref{fig6}B. 

\begin{figure}[ht]
\begin{center}
	\includegraphics[width=100mm]{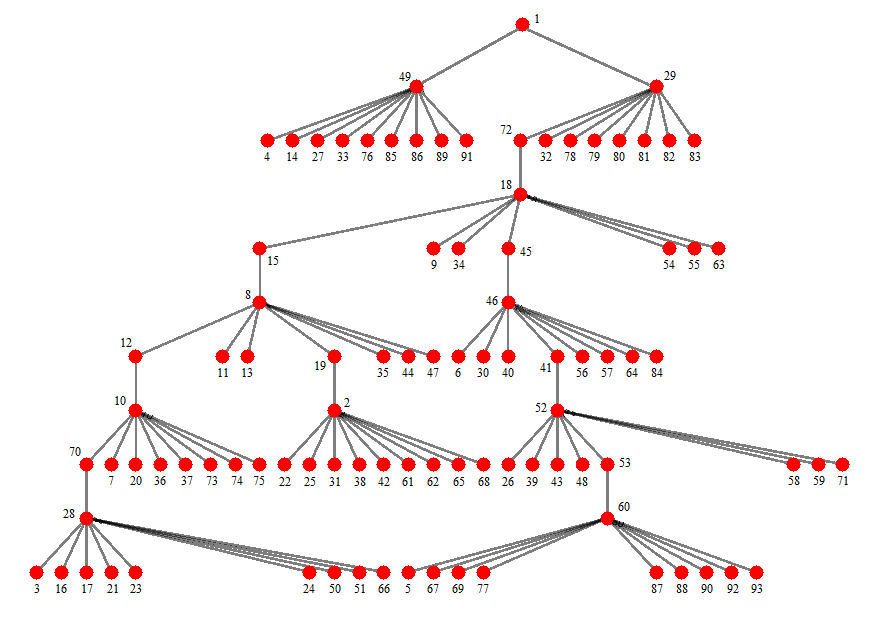}
\end{center}
\caption{{\bf Sparse spanning tree of the given gene network.}}
\label{fig7}
\end{figure}

From Fig~\ref{fig6}, it is reasonable to predict the key roles played by CUL7 and SIRT7 in this particular network as the corresponding nodes are connected to most other vertices. On the other hand, one may also argue that gene ACLY plays the most important role (in this network) as it connects the two aforementioned genes.  Similarly, the sparse spanning tree in Fig~\ref{fig7} is produced with objective value 109258. We will discuss a more systematic analysis of this network in a later section.

\subsubsection{A brain network}
\begin{figure}[htbp]
\begin{center}
	\includegraphics[width=75mm]{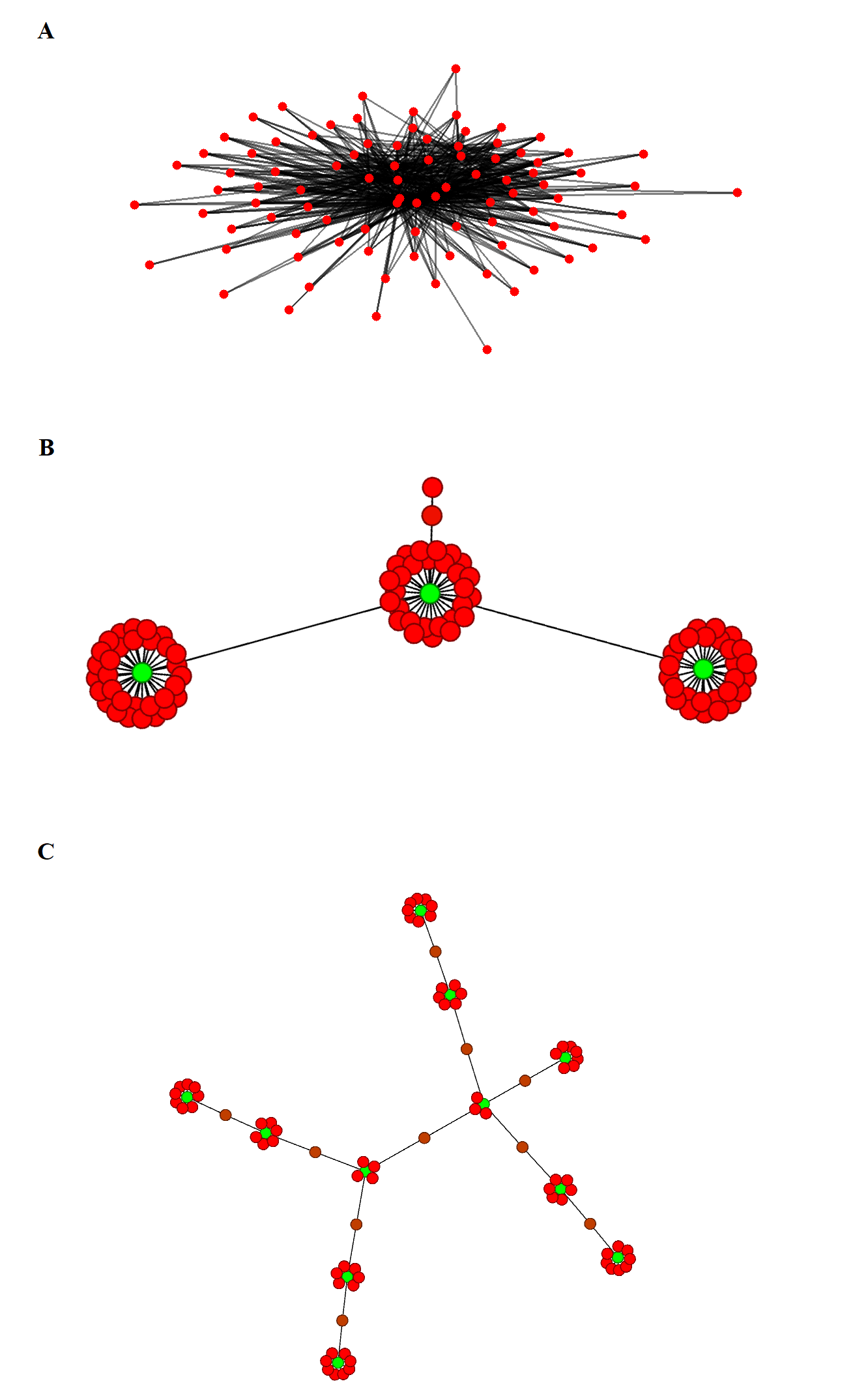}
\end{center}
\caption{{\bf The brain network figures.}
(A) The original graph. (B) The obtained dense spanning tree. (C) The obtained sparse spanning tree. Nodes of relatively high degree are highlighted in green.}
\label{fig8}
\end{figure}

Next, we consider a brain network data with 91 nodes and 628 edges \cite{22.} where nodes represent parts of the cortex and the edges represent connections between them. The original network graph is given in Fig~\ref{fig8}A. Model 2 was used with the objective function $C_{\vec{j}}$ for $\vec{j} = \langle 4,2,2,2\rangle$. We obtain the dense spanning tree in Fig~\ref{fig8}B and the sparse spanning tree in Fig~\ref{fig8}C.

If  the high-degree nodes are highlighted in the resulted spanning trees, it is easy to identify the most central parts of the cortex. This analysis can also be easily employed to study social networks as shown in the next two applications.

\subsubsection{Social networks}
\begin{figure}[htbp]
\begin{center}
	\includegraphics[width=80mm]{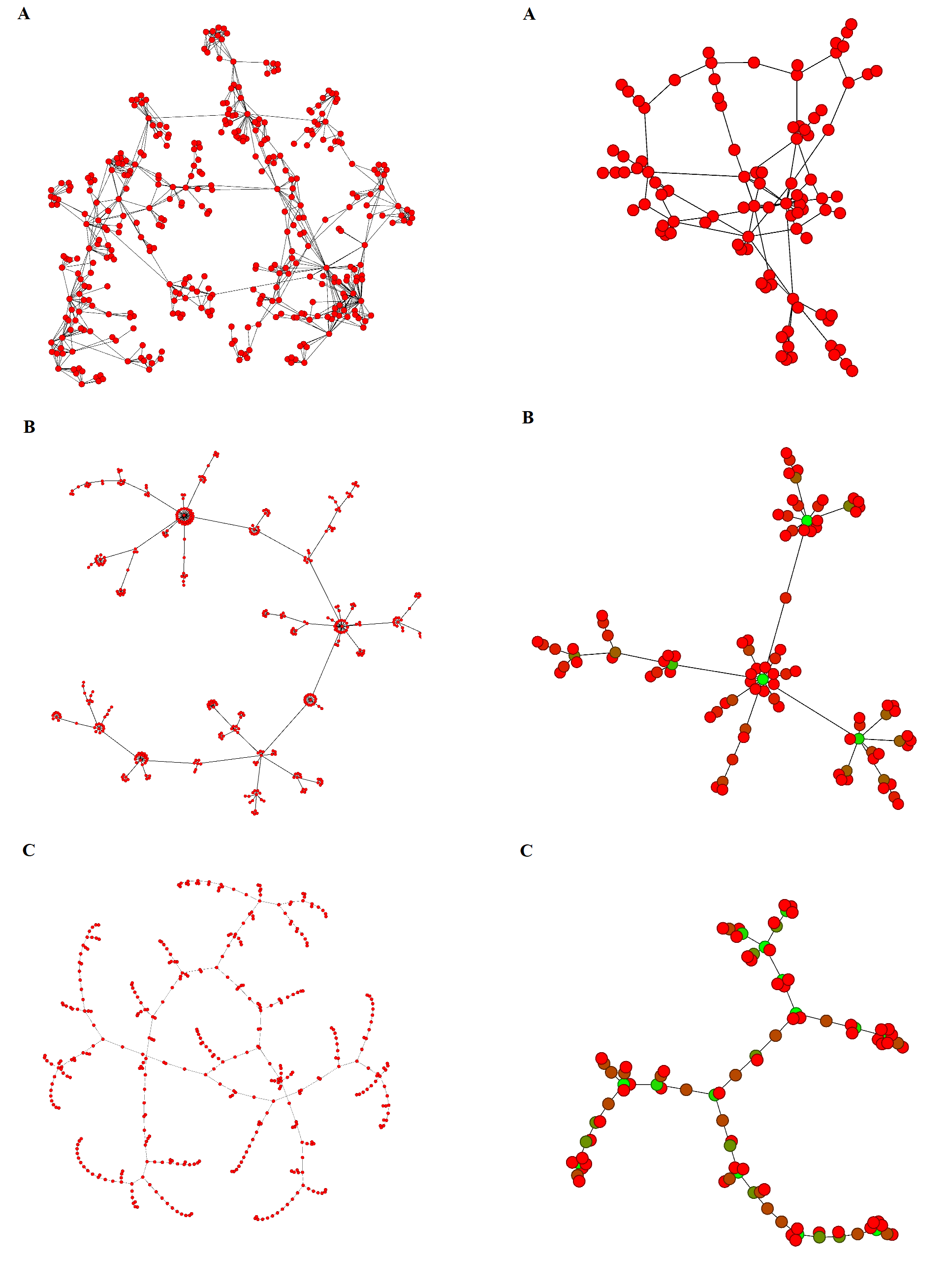} 
\end{center}
\caption{{\bf The collaboration network (left) and the re-tweet network (right).}
(A) The original graph. (B) The obtained dense spanning tree. (C) The obtained sparse spanning tree.}
\label{fig9}
\end{figure}

When applied to a collaboration network with 379 nodes and 914 edges \cite{22.,23.} and a re-tweet network with 96 nodes and 117 edges \cite{22.,24.} in Fig~\ref{fig9}, our method once again produces the dense and sparse spanning trees very quickly. As in previous example, if the high-degree nodes are highlighted, it is easy to identify the centers of these social networks.

\subsubsection{Road network}
\begin{figure}[htbp]
\begin{center}
	\includegraphics[width=50mm]{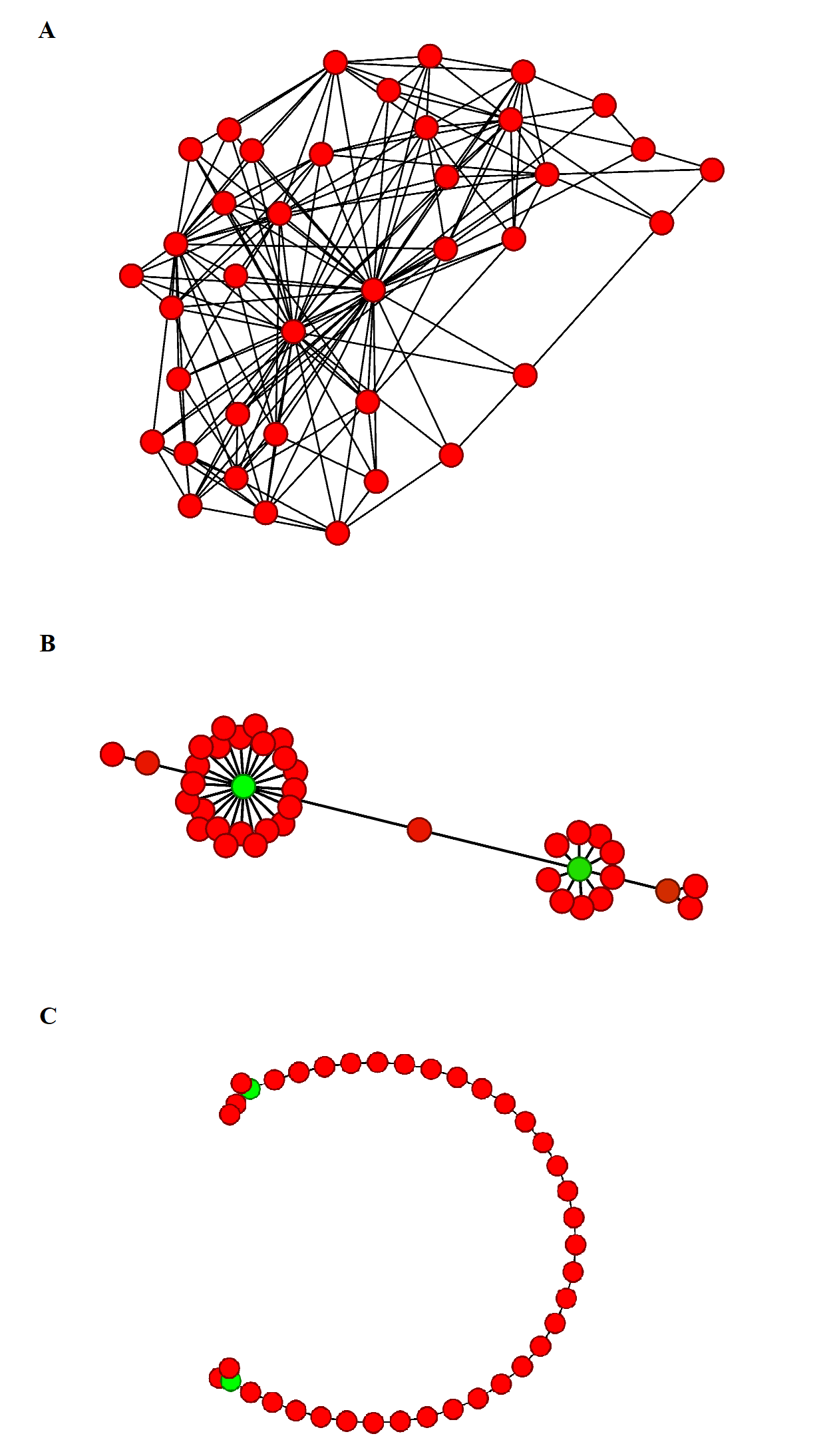}
\end{center}
\caption{{\bf The Chesapeake road network figures.}
(A) The original graph. (B) The obtained dense spanning tree. (C) The obtained sparse spanning tree.}
\label{fig10}
\end{figure}

For the last application of our method, we consider the road network of Chesapeake with 39 nodes and 170 edges \cite{22.} in Fig~\ref{fig10}A. The nodes in this network represent some locations in Chesapeake area while the edges represent the roads in between. The results are self-explanatory, identifying the centers of traffic in the dense spanning tree in Fig~\ref{fig10}B.

A summary of the above case studies is provided in Table~\ref{table2}. A comparison between the performance of Model 1 and Model 2 on the cases above could not be performed since the graphs are relatively large and therefore, only Model 2 is applied.  

\begin{table}[htbp]
\centering
\caption{
{\bf The objective function values of the applications.}}
\begin{tabular}{|c|c|c|}
\hline
& \bf Dense spanning tree & \bf Sparse spanning tree \\ 
\hline
\bf Applications & \bf Objective value &  \bf Objective value\\ \hline
Gene Network & 3675370 & 109258 \\ \hline
Brain Network & 2932444 & 80012 \\ \hline
Collaboration Network & 3266670 & 71764 \\ \hline
Re-tweet Network & 206158 & 3088 \\ \hline
Chesapeake Road Network & 372688 & 2438 \\ \hline
\end{tabular}
\label{table2}
\end{table}

\section{Discussion}
In the previous section we have applied the new GA based approach to several practical problems. It is important to note that this method is also applicable to directed graphs and graphs with certain constraints on each node for incoming and outgoing edges.

In what follows we briefly describe how  the new  approach, after minor modifications, can be used to solve  similar problems to the ones discussed in the previous section.
\begin{itemize}

\item The $k$-DST: finding the dense or sparse subtree that contains exactly $k$ vertices. 

Take, for instance, the random graph in Fig~\ref{fig3}A,  and assume $k=6$. Then, this problem can be solved using two-stage GA with sufficient constraints. The first stage is to find connected subgraphs with 6 nodes. 

In the second stage, we directly apply Model 1 or Model 2 to the subgraphs found in the first stage using the chosen objective function to find the good approximation of the globally optimal dense spanning tree. The 6-vertex dense subtree in Fig~\ref{fig11} is again obtained using Model 1 or Model 2   by maximizing the  simple objective function $S_{p=2}(\vec{h})$. 

\begin{figure}[ht]
\begin{center}
	\includegraphics[width=100mm]{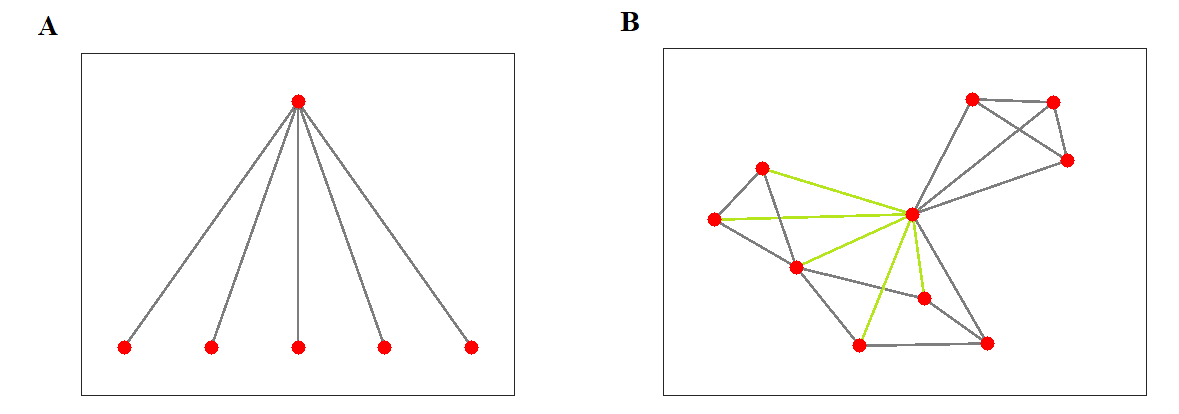}
\end{center}
\caption{{\bf Results from solving the 6-DST from the random graph in Fig~\ref{fig3}.}
(A) The optimal solution which is a star constructed by the 6 nodes of the graph. (B) The optimal solution colored in the graph, with $\vec{h} = \langle 2, 5, 8, 9, 12\rangle$.}
\label{fig11}
\end{figure}

\item The Steiner DST: finding the dense subtree that connects a given set of terminals. 

 A similar strategy used for $k$-DST problem can be applied in this case. We simply need to add constraints for the terminal nodes so that they are connected in the resulting subgraphs. Then, Model 1  or Model 2 can  be applied directly.

\item The DST with conflict pairs: given a collection of conflicting pairs of edges, find  dense subtrees which can contain at most one of the edges from each pair.

To illustrate this case, suppose that the edges labeled by 13 and 16 in Fig~\ref{fig3}A are listed as the only conflicting pair. We just need to assign very high objective values to the subgraphs containing both of these edges so that they cannot be in the optimal solution. For this example, we have the approximate optimal solution $\vec{h} = \langle 2, 5, 8, 9, 11, 12, 14, 15, 16\rangle$ as shown in Fig~\ref{fig12}. Another one  is $\vec{h} = \langle 2, 5, 8, 9, 11, 12, 13, 14, 15 \rangle$.

\begin{figure}[ht]
\begin{center}
	\includegraphics[width=100mm]{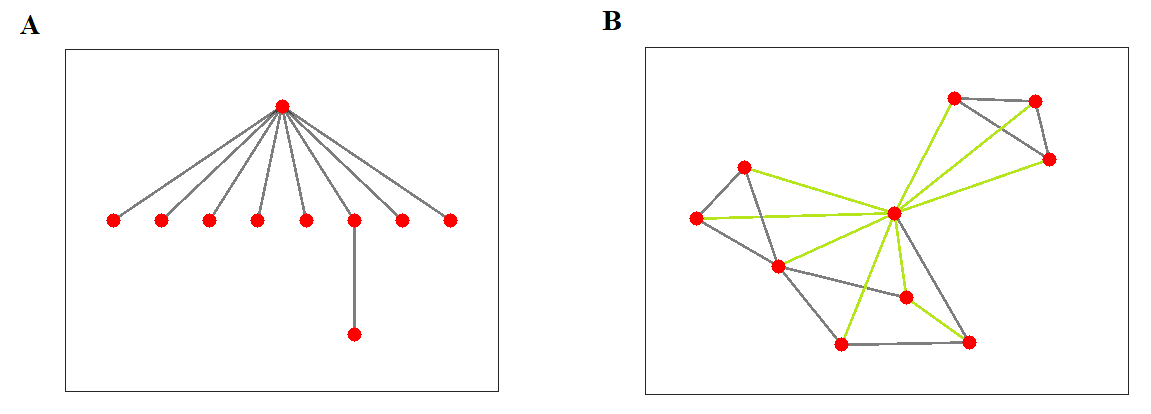}
\end{center}
\caption{{\bf DST with conflicting pairs. Applied to the random graph in Fig~\ref{fig3} where edges 13 and 16 are conflicting pairs.}
(A) The solution containing the edge 16. (B) The solution colored in the graph with $\vec{h} = \langle 2, 5, 8, 9, 11, 12, 14, 15, 16\rangle$.}
\label{fig12}
\end{figure}

\item Degree constraint DST: finding the densest spanning tree where the maximum vertex degree is bounded by a certain constant $k$. 

We only need to add the constraint of the maximum degree to the existing set of constraints and apply Model 1 or Model 2. 
\end{itemize}

As we have seen in the multiple examples above, the resulting  dense or sparse trees obtained by the new approach   may  not be  unique. To truly understand a structure, we may recursively apply our algorithm. For example, take the gene network in Fig~\ref{fig5}. The practical question of interest to biologists is usually to identify the key gene or link in the network. The dense spanning tree shown in Fig~\ref{fig6} seems to suggest the nodes with labels ``2'' and ``49'', and perhaps the node connecting them in the spanning tree, labeled ``69'' are the key genes. We now further examine the structure through the following process:
\begin{itemize}
\item[(a)] In the resulting dense spanning tree, identify the two nodes with the highest degree and remove the path connecting them from the original network;
\item[(b)] Stop if the new network is disconnected; otherwise, find the dense spanning tree from the new network;
\item[(c)] Repeat steps (a) and (b) as needed.
\end{itemize}
Table~\ref{tab:keynode} summarizes the result of this process with all removed high-degree nodes recorded. It is easy to see from the data that node ``10'', corresponding to gene HDAC5, is a third key gene that would not have been obvious from the first dense spanning tree in Fig~\ref{fig6}.

\begin{table}[htbp]
\centering
\caption{
{\bf Finding key nodes in gene networks.} For each pair $(x, y)$ in column 2, $x$ indicates node label and $y$ is the degree of $x$ in the current network. In case of three nodes with the same highest degree, two are chosen randomly before the path in between is removed.}
\scalebox{1}{
\begin{tabular}{|c|c|c|}
\hline
\bf Number of runs& \bf (Node, Degree) & \bf Removed Path \\ 
\hline
1. & (2, 42), (49, 23) & None (dense network in Fig~\ref{fig8}A) \\ \hline
2. & (2, 41), (49, 23) & 2 - 69 - 49 \\ \hline
3. & (2, 40), (29, 15), (10, 13) & 2 - 3 - 49 \\ \hline
4. & (2, 39), (49, 21) & 2 - 71 - 29 \\ \hline
5. & (2, 38), (49, 18) & 2 - 12 - 49 \\ \hline
6. & (2, 37), (49, 17), (10, 15) & 2 - 4 - 49 \\ \hline
7. & (2, 36), (49, 16), (10, 15) & 2 - 14 - 49 \\ \hline
8. & (2, 35), (29, 15), (10, 15) & 2 - 16 - 49 \\ \hline
9. & (2, 34), (29, 15), (10, 15) & 2 - 17 - 29 \\ \hline
10. & (2, 33), (49, 13), (10, 15) & 2 - 19 - 29 \\ \hline
11. & (2, 32), (29, 16), (10, 15) & 2 - 25 - 10 \\ \hline
12. & (2, 31), (49, 18), (10, 16) & 2 - 31 - 29 \\ \hline
13. & (2, 30), (49, 17), (10, 17) & 2 - 40 - 49 \\ \hline
14. & (2, 29), (49, 17), (10, 17) & 2 - 55 - 49 \\ \hline
15. & (2, 28), (49, 17), (10, 18) & 2 - 26 - 49 \\ \hline
16. & (2, 27), (49, 16), (10, 18) & 2 - 13 - 10 \\ \hline
17. & (28, 47), (49, 20) & 2 - 70 - 10 \\ \hline
18. & (28, 46), (49, 15) & 28 - 9 - 49 \\ \hline
19. & (10, 28), (49, 16), (2, 15) & 28 - 21 - 49 \\ \hline
20. & Network is disconnected & 10 - 76 - 49 \\ \hline
\end{tabular}}
\label{tab:keynode}
\end{table}

Similar process can be carried out by continuing to remove edges with the highest degree sum from their incident vertices (in the resulting dense spanning tree) from the network until the network is disconnected. The results are presented in the Table~\ref{tab:keylink}.

\begin{table}[htbp]
\centering
\caption{
{\bf Removing edges of high degree sum.}}
\scalebox{1}{
\begin{tabular}{|c|c|}
\hline
\bf Number of runs & \bf Removed Edges \\ 
\hline
1. & None (dense network in Fig~\ref{fig8}A) \\ \hline
2. & (2 - 69), (2 - 71), (49 - 69) \\ \hline
3. & (2 - 11), (2 - 3), (3 - 49) \\ \hline
4. & (2 - 25), (2 - 12), (12 - 49) \\ \hline
5. & (2 - 13), (2 - 4), (4 - 49) \\ \hline
6. & (2 - 19), (2 - 14), (14 - 49) \\ \hline
7. & (2 - 6), (2 - 41), (6 - 10) \\ \hline
8. & (2 - 16), (2 - 7), (7 - 10) \\ \hline
9. & (9 - 28), (28 - 69), (49 - 76) \\ \hline
10. & (10 - 12), (10 - 26), (26 - 49) \\ \hline
11. & (2 - 17), (2 - 33), (17 - 28) \\ \hline
12. & (28 - 47), (28 - 71), (71 - 29) \\ \hline
13. & (21 - 28), (10 - 29), (28 - 45) \\ \hline
14. & (28 - 25), (10 - 25) \\ \hline
15. & (28 - 31), (16 - 28), (16 - 49), network is disconnected \\ &  due to removal of: (16 - 28) and (16 - 49) \\ \hline
\end{tabular}}
\label{tab:keylink}
\end{table}

\newpage
\section{Conclusion}

The novel use of the GA  presented in this paper successfully employs the degree conditions as the new criteria to find dense and sparse spanning trees for any connected, directed or undirected, weighted or unweighted graphs. Compared with the previously established heuristic algorithm for the same purpose, the GA, in particular Model 2 which is enhanced by the use of Kruskal algorithm is more efficient and obtains approximate solutions that are much closer to the optimal solutions. In addition to outperforming former algorithms, the proposed approach can also be easily adapted to solve similar problems under various additional constraints.

We also  discuss novel recursive application of  the new method to obtain a deeper understanding of the graph structure. Using gene networks as an example, our approach finds key genes in a network that were not obvious from previous studies.

We conclude the paper with an important observation: in the above mentioned approaches, the number of runs it takes to disconnect the network (as recorded in Tables~\ref{tab:keynode} and \ref{tab:keylink}) also measures how strongly the network is connected. Hence, the ``density'' of a network can be measured from this novel perspective, which is more accurate than  some other  simple criteria such as the number of edges.

\section*{Acknowledgement}
The authors declare that no conflict of interests exist.


\begin{thebibliography}{}


\bibitem{1.}
Gargano, L., Hell, P., Stacho, L., Vaccaro, U.: Spanning trees with bounded number of branch vertices.In 29th International Colloquium on Automata, Languages and Programming (ICALP). Lecture Notes in Computer Science, Springer, Berlin, 2380, 355-365 (2002) 

\bibitem{2.}
Silva, R., Silva, D., Resende, M., Mateus, G., Goncalves, J., Festa, P.: An edge-swap heuristic for generating spanning trees with minimum number of branch vertices. Optim. Lett. 8, 1225-1243 (2014)  

\bibitem{3.}
Bazlamacci, C., Hindi, K.: Minimum-weight spanning tree algorithms: A survey and empirical study. Computers and Operations Research 28, 767-785 (2001)

\bibitem{4.}
Hwang, F. K., Richards, D. S., Winter, P.: The steiner tree problem. North-Holland, New York (1992)

\bibitem{5.}
Narula, S. C., Ho, C. A.: Degree-constrained minimum spanning tree. Comput. Oper. Res. 7, 239-249 (1980)

\bibitem{6.}
Ozen, M., Wang, H., Wang, K., Yalman, D.: An edge-swap heuristic for finding dense spanning trees. Theory Appl. Graphs. 3(1), 1-10 (2016)

\bibitem{7.}
Darmann, A., Pferschy, U., Schauer, J.: Minimal spanning trees with conflict graphs. Optimization online (2009)

\bibitem{8.}
Amberg, A., Domschke, W., Voss, S., Vo, S.: Capacitated minimum spanning trees: algorithms using intelligent search. Combinatorial Optimization: Theory and Practice 1, 9-40 (1996)

\bibitem{9.}
Li, T., Gao, Y., Dong, Q., Wang, H.: Degree sums and dense spanning trees. PLoS ONE, 12(9): e0184912. https://doi.org/10.1371/journal.pone.0184912 (2017)

\bibitem{10.}
Wiener, H.: Structural determination of paraffin boiling points. J. Am. Chem. Soc. 69, 17-20 (1947)

\bibitem{11.}
Wiener, H.: Correlation of heats of isomerization, and differences in heats of vaporization of isomers among the paraffin hydrocarbons. J. Am. Chem. Soc. 69, 2636-2638 (1947)

\bibitem{12.}   Sz\'ekely, L. A., Wang, H.: On subtrees of trees. Adv. Appl. Math. 34, 138-155 (2005) 

\bibitem{13.}  Wagner, S.: Correlation of graph-theoretical indices. SIAM J. Discrete Mathematics 21(1), 33-46 (2007)

\bibitem{14.} Schmuck, N., Wagner, S., Wang, H.: Greedy trees, caterpillars, and Wiener-type graph invariants. MATCH Commun.Math.Comput.Chem. 68(1), 273-292 (2012)  
\bibitem{15.} Wang, H.: The extremal values of the Wiener index of a tree with given degree sequence, Discrete Applied Mathematics 156, 2647-2654 (2008)

\bibitem{16.} Zhang, X. D., Xiang, Q. Y., Xu, L.Q., Pan, E. Y.: The Wiener index of trees with given degree sequences, MATCH Commun.Math.Comput.Chem. 60, 623-644 (2008)

\bibitem{17.} Zhang, X. M., Zhang, X. D., Gray, D., Wang, H.: The number of subtrees of trees with given degree sequence, J. Graph Theory 73(3), 280-295 (2013)

\bibitem{18.}
Mitchell, M.: An introduction to genetic algorithms. MIT Press, Cambridge, MA (1996)

\bibitem{19.}
MathWorks: Global optimization toolbox: User's guide (R2018b). Retrieved Sept 17, 2018 from \url{https://www.mathworks.com/help/pdf_doc/gads/gads_tb.pdf}

\bibitem{20.}
Kruskal, J. B.: On the shortest spanning subtree of a graph and the traveling salesman problem. Proceedings of the American Mathematical Society 7, 48-50 (1956)	

\bibitem{21.}
Cormen, T., Leiserson, C. E., Rivest, R. L., Stein, C.: Introduction To Algorithms. MIT Press (3rd ed.), ISBN 0262258102 (2009)

\bibitem{25.}
IBM ILOG CPLEX Optimization Studio, Retrieved Jan 18, 2020 from \url{https://www.ibm.com/products/ilog-cplex-optimization-studio}

\bibitem{22.}
Rossi, A. R., Ahmed, N. K.: The network data repository with interactive graph analytics and visualization. Proceedings of the Twenty-Ninth AAAI Conference on Artificial Intelligence, (2015) \url{http://networkrepository.com}

\bibitem{23.}
Newman, M. E. J.: Finding community structure in networks using the eigenvectors of matrices. Physical review E 74(3) (2006)

\bibitem{24.}
Rossi, R. A., Gleich, D. F., Gebremedhin, A. H., Patwary, M. A.: What if clique were fast? Maximum cliques in information networks and strong components in temporal networks.arXiv preprint arXiv:1210.5802, 1-11 (2012)

\end{thebibliography}
\end{document}